\documentclass{article}
 \usepackage{amsmath,amsthm} 
     \usepackage{amssymb} 
     \usepackage{enumerate} 
     \usepackage{bm} 
     \usepackage{hyperref} 
\input xy
\xyoption{all}

 \newcommand{\PP}{\mathbf{P}} 
   \newcommand{\ZZ}{\mathbb{Z}} 
    
   \newcommand{\OO}{\mathcal{O}} 
    
      \newcommand{\OP}{\mathcal{O}_{\PP^2}} 
   \newcommand{\OPN}{\mathcal{O}_{\PP^N}}

   \newcommand{\II}{\mathcal{I}} 
   \newcommand{\FF}{\mathcal{F}} 
   \newcommand{\EE}{\mathcal{E}} 
   \newcommand{\TT}{\mathcal{T}} 
   \newcommand{\GG}{\mathcal{G}} 
   \newcommand{\HH}{\mathcal{H}}

   \newcommand{\IY}{\II_{Y}}

   \newcommand{\gsr}{\mathrm{gsrk}} 
      \newcommand{\coker}{\mathrm{coker}} 
   \newcommand{\rk}{\mathrm{rk}}

\newcommand{\im}{\mathrm{Im}\kern.3pt} 
\newcommand{\st}{\mathrm{st}} 
\newcommand{\gst}{\mathrm{gst}}

    \newtheorem{definition}{Definition}
    \newtheorem{lemma}{Lemma}
    \newtheorem{proposition}{Proposition}
    \newtheorem{theorem}{Theorem}
    \newtheorem{corollary}{Corollary}
    \newtheorem{example}{Example}

\title{Splitting type, global sections  and Chern classes for  torsion free sheaves on $\PP^N$} 

   \author{ Cristina Bertone-Margherita Roggero \\ 
                   } 

\date{\empty}

\begin{document}
\maketitle

\begin{abstract}
In this paper we compare a  torsion free sheaf $\FF$ on $\PP^N$ and the free vector bundle $\oplus_{i=1}^n\OPN(b_i)$ having same rank and  splitting type. We show that the first one has always \lq\lq less\rq\rq\ global sections, while it has a higher second Chern class. In both cases bounds for the difference are found in terms of the maximal free subsheaves of $\FF$.   As a consequence we obtain a direct, easy and more general proof of the \lq\lq Horrocks' splitting criterion\rq\rq,  also holding for torsion free sheaves, and lower bounds for the  Chern classes $c_i(\FF(t))$ of twists of $\FF$,  only depending on some numerical invariants of $\FF$. Especially, we prove  for rank $n$ torsion free sheaves on $\PP^N$, whose splitting type has no gap (i.e.  $b_i\geq b_{i+1}\geq b_i-1$ for every $i=1, \dots,n-1$ ),    the following formula for the discriminant:
\[ \Delta(\FF):=2nc_2-(n-1)c_1^2\geq -\frac{1}{12}n^2(n^2-1)\]  
Finally in the case of rank $n$ reflexive sheaves  we obtain polynomial upper bounds for the absolute value of the higher Chern classes $c_3(\FF(t)), \dots, c_n(\FF(t))$,  for the dimension of the cohomology modules $H^i\FF(t)$ and for the Castelnuovo-Mumford regularity of $\FF$; these polynomial bounds  only depend only on $c_1(\FF)$, $c_2(\FF)$,  the splitting type of $\FF$ and $t$.
\end{abstract}



\section*{Introduction}

The present paper deals with the  problem of determining  all the possible values for the Chern classes of a rank $n$ torsion-free sheaf $\FF$ on a projective space $\PP^N$. To this aim
 we consider in particular  totally split  bundles which are  in some sense  close to $\FF$ and  find relations between the Chern classes of  $\FF$ and those, very easy to compute,  of the split bundles.  In the previous paper \cite{prec} we studied  maximal free subsheaves  of a  reflexive sheaf $\FF$; here we are mainly concerned with  the sheaf $\oplus_{i=1}^n\OPN(b_i)$ ($\OPN(\mathbf{b})$ for short), whose restriction to a general line $L$ is   isomorphic to $\FF_L$. 

This method, though rather elementary, is in fact direct and effective and allows very general results. 
 As an evidence of the efficiency of our method, we point out for instance
  the simple, self-included   and constructive proof of the  splitting criterion for torsion free sheaves  given in \S 3 (Theorem \ref{previous}) which generalizes \lq\lq Horrocks'  criterion\rq\rq\ for vector bundles:
 
 \medskip

 \noindent \textbf{ Theorem A.} \textit{Let  $\FF$ be a torsion free sheaf  on $\PP^N$, $N\geq 2$, $H$ be a hyperplane in $\PP^N$ and $\mathbf{b}=[b_1, \dots,b_n]$ be an ordered sequence of integers.}

\textit{Then  $\FF\simeq \OPN (\mathbf{b})$ if and only if $\FF_H\simeq \OO_H(\mathbf{b})$ and either one of the following conditions holds:}

\textit{	i)   $N=2$   and  $H^1 (\FF(t))=0$ for every $t\leq -b_n$}
	
\textit{	ii) $N\geq 3$    and  $H^1 (\FF(t))=0$ for every $t\ll 0$. }

 \medskip
 
 \noindent (for a different proof of a similar result see \cite{ABE}).
The   sequence of integers $\mathbf{b}=[b_1,\dots, b_n]$ is the so-called  \textit{splitting type} of $\FF$,     a classical invariant in both  algebraic  and  differential theory of vector bundles,  which is directly connected to important properties like uniformity and stability  (see for instance \cite{okverde}). Moreover $\mathbf{b}$  determines the first Chern class of $\FF$; it is well known that    $c_1(\FF)$ is simply the sum of the integers $b_i$, that is $c_1(\FF)=\sum b_i=c_1(\OPN(\mathbf{b}))$. Very little is known about relations between the splitting type and the other Chern classes. One could think that no such a relation exists, because the splitting type can be defined only   restricting  $\FF$ to a  line, where the higher Chern classes  $c_s$,  $s\geq 2$, disappear. 

However some interesting inequalities can be found in literature for both vector bundles of any rank on $\PP^n$ and rank $2$ reflexive sheaves on $\PP^3$ (see \cite{EF}, \cite{stable} and \cite{sawer}) that  we generalize (Corollary \ref{risH0}, Theorem \ref{secondobound}, and  Corollary \ref{alte}):  

\medskip

\noindent \textbf{Theorem B.} \textit{Let $\FF$ be a   rank $n$ torsion free sheaf on $\PP^N$. Then:}
\begin{itemize}
	\item[\textit{i)}] $h^0\OPN(\mathbf{b})-h^0\FF\geq 0$  
	\item[\textit{ii)}] 	$c_2(\FF)\geq c_2(\OPN(\mathbf{b}))=\sum b_ib_j$
	\item[\textit{iii)}] $c_s(\FF(t))\geq c_s(\OPN(\mathbf{b}\otimes \OPN(t))$ \textit{for every }$s\leq n$ and $t\gg 0$ 
\end{itemize}

\medskip 
 
 Actually what we prove is something more precise than the simple positivity of the three differences $h^0\OPN(\mathbf{b})-h^0\FF$, $c_2(\FF)- c_2(\OPN(\mathbf{b}))$ and $c_s(\FF(t))-c_s(\OPN(\mathbf{b}+t))$; in fact we obtain   lower bounds for them that involve the maximal free subsheaves of $\FF$ and of its restriction $\FF_H$ to general linear subspaces $H$ in $\PP^N$. Every free subsheaf of $\FF$ is, in some sense, also a subsheaf of $\OPN(\mathbf{b})$, while the converse can hold only if $\FF$ itself is totally split (Corollary \ref{primospezz}; see also  \cite{prec}). As a consequence  we show that equality can realize  in either one of \textit{i)}, \textit{ii)}, \textit{iii)}  if and only if $\FF\simeq \OPN(\mathbf{b})$.

It is not difficult to see that  for every   $  s\leq n$ the difference  $\delta_s(t)=c_s(\FF(t))-c_s(\OPN(\mathbf{b}+t))$ can be expressed as a polynomial in $t$ and coefficients depending on the  Chern classes $c_1, \dots, c_{s}$, whose degree with respect to $t$ is at most $s-2$;  the fact that equality cannot hold in \textit{ii)} unless $\FF\simeq \OPN(\mathbf{b})$, allows us to conclude    that for every torsion free sheaf which is not totally split the $t$-degree of $\delta_s(t)$ is precisely $s-2$; moreover    our lower bound for $\delta_s(t)$  is a degree $s-2$ polynomial too. An analogous property holds with respect to the difference at the left hand of \textit{i)} and to the relative lower bound: they are polynomials of the same degree    $N-2$.

 Interesting consequences of the above quoted results are the following  lower bound on the 
 discriminant $\Delta (\FF)=2nc_2-(n-1)c_1^2$ of a   rank $n$ torsion free sheaf $\FF$ on $\PP^N$ (Theorem \ref{compleanno}):
 
 \medskip
 
\noindent \textbf{Theorem C. } \textit{ If $\FF$   is a torsion free sheaf and its splitting type has no gap,  that is  $b_i\geq b_{i+1}\geq b_i-1$ for every $i=1, \dots,n-1$,  then} $ 12\Delta (\FF)\geq - n^2(n^2-1).$

\textit{Moreover $  \Delta (\FF)\geq  2n$ if a plane section of $\FF$ is     semistable and $\Delta (\FF)\geq \frac34 n^2$
 if the plane section is stable.}

	\medskip

	When $n=2$ this formula gives Schwarzenberger inequality for rank 2 semistable vector bundles and,  for $N=2$ and higher $n$,  it is close to   Bogomolov inequality  $\Delta (\FF)\geq 0$  for semistable sheaves  (\cite{bog}),  generalized to every $N$  using  semistability results  of the restriction to suitable surfaces (\cite{gie}, \cite{HL}). 
	
	\medskip

Finally we generalize to reflexive sheaves  the upper-bounds for vector bundles obtained in \cite{EF} (see  Theorem 3.3 and Theorem 4.2);they concern the dimension of the cohomology  modules $H^i\FF$, the absolute value of the Chern classes $c_s(\FF)$ and the regularity  of $\FF$ and they are obtained through   polynomial functions depending on the rank of $\FF$, on the dimension of the projective space, on the first and second Chern classes $c_1(\FF)$ and $c_2(\FF)$ and on the splitting type $\mathbf b$ (see Theorem \ref{polinomi}), while no such bounding formulas  exist for the wider class of torsion free sheaves.
 
 \medskip

\noindent \textbf{Theorem D.}
 \textit{For any choice of non-negative integers $n$, $N$ and $s$, with $N\geq 2$ and $3\leq s \leq N$, there are suitable polynomial fuctions $P_{n,N}$, $Q_{n,N}$, $C_{n,N,s}$  in the set of variables $\{ c_1,c_2, \mathbf{b}, d, \delta_2\}$ (or either one of the sets of Lemma \ref{sottoinsieme}), such that for any reflexive sheaf $\FF$ of rank $n$ on $\PP^N$ we have:}
 \begin{enumerate}
	\item[$i)$] \textit{$h^i\FF \leq P_{n,N}$ for all $0\leq i\leq N$;
	\item[ii)]  $h^i\FF(k)=h^1\FF(-k)=h^0\FF(-k)=0$ for all   $k \geq Q_{n,N}$ and $1\leq i\leq N$;
	\item[iii)] the Castelnuovo-Mumford regularity of $\FF$ is lower than  $ Q_{n,N}$ and especially $\FF(k)$ is generated by global sections for $k \geq Q_{n,N}+N$;
	\item[iv)] $\vert c_s(\FF)\vert  \leq C_{n,N,s}$.}
\end{enumerate}

 \medskip

 In \S 1 we collect some basic definitions and  known results, that we will use more often in the following. 
\S 2, \S 4 and  \S 6 concern respectively  \textbf{Theorem B} \textit{i), ii)} and \textit{iii)}.  
In \S 3 we prove Horrocks' criterion for torsion free sheaves, that is  \textbf{Theorem A},  in \S 5  the lower bounds for the discriminant that is  \textbf{Theorem C} and in \S 7 the polynomial upper bounds for cohomology, Castelnuovo-Mumford regularity and higher Chern classes that is  \textbf{Theorem D}.  
Finally in \S 8 we exhibit  examples showing the sharpness of the main results of the paper.

\section{Generalities}

We consider an algebraically closed field $k$ of characteristic 0;  
$\PP^N$ is the pro\-jec\-ti\-ve space of dimension $N$ over $k$. As usual, if $\FF$ is a coherent sheaf on $\PP^N$, we will denote by $h^i(\FF)$  the dimension of the $i$-th cohomology module $H^i(\FF)$ as a $k$-vector space and by $H^i_*\FF$ the direct sum $\oplus_{n\in \ZZ} H^i\FF(n)$; in particular $H^0_*\OPN=k[X_0,\dots,X_n]$ and, in a natural way,  $H^0_*\FF$ is a $H^0_*\OPN$-module. 
 $\widehat{\FF}$ will denote the subsheaf of a sheaf $\FF$ generated by  ${H^0\FF}$. 

\medskip

\noindent \textbf{1)}  For every   rank $n$  coherent sheaf $\FF$  on $\PP^N$, we denote by $c_i(\FF)$ or 
simply $c_i$ ($i=1, \dots , N$) its Chern classes  that we think  as integers 
and by 
$
C_t(\FF)=\sum_{i=0}^N c_i(\FF)t^i
$ its Chern polynomial. 
If 
$0\rightarrow \FF'\rightarrow \FF \rightarrow \FF''\rightarrow 0$ is an exact sequence, then 
$ 
C_t(\FF)=C_t(\FF')C_t(\FF'') 
$ in $  \ZZ[t]/(t^{N+1})$. 
Moreover, for every $l\in \ZZ$, the Chern classes 
of $\FF (l)$ are given by: 
\begin{equation}\label{cherncontwist} 
c_i(\FF(l))=c_i+(n-i+1)lc_{i-1}+{n-i+2\choose 2}l^2c_{i-2}+\cdots+{n\choose 
i}l^i. 
\end{equation} 
The Chern classes of a rank $n$ split bundle $\oplus_{i=1}^n\OPN(b_i)$ are $c_l=\sum b_{i_1}b_{i_2}\cdots b_{i_l}$ (where the sum is over all  $i_1 <i_2<\dots <i_l$), if $l\leq min\{N,n\}$ and $c_l=0$ otherwise.

 The discriminant of the rank $n$ torsion free sheaf $\FF$ on $\PP^n$ is  $\Delta (\FF)=2nc_2-(n-1)c_1$ (for us  integer numbers).
 
 If $H$ is a general hyperplane, the torsion free sheaf $\FF$ and its  restriction $\FF_H$ to $H$ are connected by the standard exact sequence:
 \begin{equation}\label{restrizionestandard} 0 \to \FF(-1) \to \FF \to \FF_H \to 0. \end{equation}
 $\FF_H$ is torsion free too and $c_i(\FF_H)=c_i(\FF)$ for every $i\leq N-1$. 
 \medskip
 
 \noindent \textbf{2)}\label{reflexive} A coherent sheaf $\FF$ on $\PP^N$  is \emph{reflexive} if 
the canonical morphism $\FF\rightarrow\FF^{\vee\vee}$ is an isomorphism, where $\FF^{\vee}$ is the dual sheaf, that is 
$\FF^{\vee}=\mathcal{H}om(\FF, \OPN )$. Vector bundles and the dual of any coherent sheaf are reflexive; for every 
integer $l$,  $\OPN (l)$ is the  the only rank 1 reflexive sheaf on $\PP^N$ 
with $c_1=l$.  We refer to \cite{stable} for general facts about reflexives sheaves. Now we only recall that for a reflexive sheaf  $\FF$ on $\PP^N$, $N\ge   2$, the first cohomology module $H^1_*\FF(t)$ is a vector space of finite dimensions because $H^1\FF(t)=0$  for $t\ll 0$: this can be deduced for instance from Corollary 1.5 of \cite{stable}.

\medskip

\noindent \textbf{3)}\label{minore}  We will denote by $\mathbf{a}=[a_1,\dots,a_n]$ a  sequences of integers such that $a_1 \geq  \dots \geq a_n$ and use the compact notation $\OPN( \mathbf{a})$ for the split bundle $\oplus_{i=1}^n \OPN(a_i)$;  moreover if  $\mathbf{a'}=[a'_1,\dots,a'_n]$ is another sequence,  we will write  $\mathbf{a} \leq \mathbf{a'}$ if    $a_i\leq a'_i$ for $i=1,\dots,n$ and $\mathbf{a} < \mathbf{a'}$  if ``$\leq$'' holds and $\mathbf{a} \neq \mathbf{a'}$, that is $a_i<a_i'$ at least once. Finally if $t\in \ZZ$, $\mathbf{a}+t=[a_1+t,\dots,a_n+t]$.
 
 \medskip
 
\noindent \textbf{4)}\label{not4} For a general hyperplane $H$, the restriction $\FF_H$ of a torsion free (reflexive) sheaf $\FF $ is torsion free (reflexive) too and the  restriction $\FF_L$ to a general line $L$ is a free vector bundle.
 We denote by $\st(\FF)$ the  {\bf 
splitting type } of $\FF$ that is the sequence of integers $\mathbf{b}=[b_1, \dots , b_n]$ ($b_i \geq b_{i+1}$) such that $\FF_L =  \OO_L (\mathbf{b})$ for a
general line $L$ in $\PP^N$; recall that 
$c_1(\FF)=b_1 + \cdots +b_n$. We say that the splitting type $\mathbf b$ has no gap if  $b_i-b_{i+1}\leq 1$ for every $i=1, \dots,n-1$.

\begin{definition}\label{typeglobalsections} 
Let $\HH$ be a coherent  sheaf on $\PP ^N$. The \textbf{rank of $\HH$ by global sections} $\gsr (\HH)$ is the   maximum $m$ for  which   there is an injective map
$
\phi :   \OPN^{m} \rightarrow \HH. 
$
The \textbf{ 
global section type} $\gst (\HH)$ is the sequence of 
integers $\mathbf{a} =[a_1, \dots, a_n]$ such that $\gsr(\HH (-a_i -1)) <i $   and 
$ \gsr(\HH (-a_i )) \geq i$ for every 
$i=1, \dots , n$. 
\end{definition} 

It is easy to see that  $a_i\geq 0$ if and only if $\gsr (\FF) \geq i$ and  $\gst(\FF (l))=[a_1+l, \dots, a_n+l]$. 
Note that  for    $\gst(\FF_L)$ for a general line $L$ is the  splitting type  $\st(F)$.  

The following Lemma collect some straightforward consequences of the above definition (for the first one see also \cite{prec} Remark 1 and Lemma 5 where this notion was first introduced).

\begin{lemma}\label{lemmafacile} Let $\FF$ be a rank $n$ torsion free  sheaf on $\PP^N$.

 \noindent $1.$ If $\widehat{\FF}$ is the subsheaf of $\FF$ generated by $H^0\FF$, then   $\gsr (\FF)=\gsr (\widehat{\FF})=\rk (\widehat{\FF})$.

\noindent $2.$ If $\mathbf{a}=\gst (\FF)$ there is a (not unique) injective map:  $\phi \colon \OPN (\mathbf{a}) \hookrightarrow \FF$ which is  \lq\lq maximal\rq\rq \  among the 
maps of that  type in the following sense: \\ 
 if $\bm{\alpha}=[\alpha_1, \dots, \alpha_r] $ and  $f \colon \OPN (\bm{\alpha})\hookrightarrow \FF$ is an injective map,  then  $r\leq n$ and $ \alpha _i \leq a_i $  for every  $i=1, \dots r$.

\noindent 3. 
 If $H \cong \PP^r$ is a  general linear space in $\PP^N$, then  $\gst 
(\FF)\leq \gst (\FF_H)$ (as a sheaf on $H$) and especially  $\gst(\FF) \leq \st (\FF)$.
\end{lemma}


\section{$h^0\FF \leq h^0\OPN(\mathbf{b})$}\label{2}

It is not difficult to see that   every torsion free sheaf  $\FF$ on $\PP^N$ has always \lq\lq less\rq\rq\ global sections than the split vector bundle $\OPN(\mathbf{b})$ having the same rank and  splitting type. 

\begin{lemma}\label{coomologiaestrema} Let $\FF$ be a rank $n$ torsion free sheaf on $\PP^N$. Then:
\begin{itemize}
	\item[1.] $h^0\FF \leq h^0\OPN(\mathbf{b})$
	\item[2.] $h^N\FF \leq h^0\OPN(-\mathbf{b}-N-1)$
\end{itemize}
\end{lemma}
\begin{proof} Both inequalities can be proved by induction on $N$ using the canonical exact sequence (\ref{restrizionestandard}). The initial case $N=1$ is obvious because $\FF = \OPN(\mathbf{b})$.
Assume $N\geq 2$ and the inequalities true on a general hyperplane $H$.  From the cohomology exact sequence of (\ref{restrizionestandard}) we obtain for every $k$:
 $$h^0\FF(k)-h^0\FF(k-1)\leq h^0\FF_H(k)\leq h^0\OO_H(\mathbf{b}+k)= h^0\OPN(\mathbf{b}+k)-h^0\OPN(\mathbf{b}+k-1).$$ 
  As $h^0\FF(k)=0$ for  $k\ll 0$, we obtain $1.$ summing up on $k\leq 0$. In the same way we obtain $2.$ summing up on $k\leq 0$ the inequalities 
 $h^N\FF(k)-h^N\FF(k+1)\leq h^{N-1}\FF_H(k+1)\leq h^0\OO_H(-\mathbf{b}-k-N-1)=h^0\OPN(-\mathbf{b}-k-N-1)-h^0\OPN(-\mathbf{b}-k-N-2)$  and recalling that $h^N\FF(k)=0$ for $k \gg 0$.
\end{proof}
In order to  find a better  lower bound for $h^0 \OPN (\mathbf{b})-h^0\FF$, we take into consideration the \lq\lq growth\rq\rq \ of the global  section type of $\FF_H$ as a function of the codimension of the general linear space $H$. 

\medskip

\begin{lemma}\label{primosezioni}
Let $\GG$  be a  subsheaf of a torsion-free sheaf $\EE$  on $\PP^N$.   Assume that $\gsr (\GG)=r $ and $\gsr (\EE)=r'$.
Then $\displaystyle h^0 \EE \geq h^0\GG + \sum_{r<i\leq r'}h^0\OPN(\alpha_i)$ where $[\alpha_1,\dots,\alpha_n]=\gst(\EE)$.
\end{lemma}
\begin{proof}
Thanks to Lemma  \ref{lemmafacile}, \textit{1.},  it will sufficient to show the inequality with respect to    $\widehat{\GG}$ and $\widehat{\EE}$ so that $\gsr(\widehat{\GG})=\rk(\widehat{\GG})=r$ and $\gsr(\widehat{\EE})=\rk(\widehat{\EE})=r'$. 

Using $\OPN^r \hookrightarrow \widehat{\GG}$,  $\widehat{\GG}\hookrightarrow \widehat{\EE}$ and $\oplus_{i=1}^{r'}\OPN(\alpha_i) \hookrightarrow \widehat{\EE}$   we can  define a map \[ \phi \colon  \oplus_{i=1}^{r'}\OPN(\alpha_i) \oplus \OPN^r \rightarrow \widehat{\EE}\]
 whose image has rank $r'$.  Now we can choose $r'-r$ summand $ \oplus_{j=1}^{r'-r}\OPN(\overline{\alpha}_j)$  in $ \oplus_{i=1}^{r'}\OPN(\alpha_i)$ such that  $\phi '\colon  \oplus_{j=1}^{r'-r}\OPN(\overline{\alpha}_j) \oplus \OPN^r \rightarrow \widehat{\EE}$  is injective and then $Im(\phi')$ has rank $r'$. Let    $\psi \colon  \oplus_{j=1}^{r'-r}\OPN(\overline{\alpha}_j) \oplus \widehat{\GG} \rightarrow \widehat{\EE}$ be the  induced map; its image, containing $Im(\phi')$, has maximal rank $r'$ and so its  kernel is a rank $0$ torsion free sheaf,  that is $\psi$   is injective.
 
 In order to conclude, we only have to observe that $[\overline{\alpha}_{1},\dots,\overline{\alpha}_{r'-r}]$ is an ordered  sublist of $[{\alpha}_{1},\dots,{\alpha}_{r}]$ and then $\overline{\alpha}_{j}\geq {\alpha}_{r+j}$.
\end{proof}

\begin{proposition}\label{utile}
Let $\FF$ be a rank $n$ torsion free sheaf on $\PP^N$ and $H$ be a general hyperplane. If    $\gsr (\FF)=r$,    $\gst(\FF_H)=[\alpha_1,\dots,\alpha_n]$, then:
\begin{equation}\label{h0Fsenzasomma}
h^0\FF-h^0\FF(-1)\leq h^0\FF_H-\sum_{ i> r}h^0\OO_H(\alpha_i).
\end{equation}
\end{proposition}

\begin{proof} Let  $(H^0\FF)_H $ be the image of the restriction map $H^0\FF \rightarrow H^0\FF_H$.  We consider the sheaf 
  $\EE=\widehat{(\FF_H)}$  and   its subsheaf $\GG$ generated  by  $(H^0\FF)_H $. Then $h^0\EE=h^0\FF_H$ and  $\rk(\EE)=\gsr(\EE)=\gsr(\FF_H)=r'\geq r$,  so that $\gst (\EE)=[\alpha_1,\dots,\alpha_{r'}]$, where $\alpha_{r'}\geq 0$ and $\alpha_j <0$ if $j>r'$; moreover  $\rk( \GG)=\gsr (\GG)=r$ and $h^0\GG \geq h^0\FF-h^0\FF(-1)$ because $\GG$  is a subsheaf of the rank $r$ sheaf $(\widehat{\FF})_H$ and $H^0\GG \supseteq (H^0\FF)_H $.

If we apply Lemma \ref{primosezioni} to the inclusion $\GG \hookrightarrow \EE$ we get: \\
$h^0\FF_H=h^0\EE \geq h^0\GG + \sum_{i=r+1}^{r'}h^0\OO_H(\alpha_i)\geq h^0\FF-h^0\FF(-1)+\sum_{i>r}h^0\OO_H(\alpha_i)$.





\end{proof}

\begin{theorem}\label{H0}
Let $\FF$ be a rank $n$ torsion free sheaf on $\PP^N$ ($N\geq 1$). For every  general linear space $H_j$ of dimension $j$ in $\PP^N$, let   $\mathbf{a}_{j}=[a_{j,1}, \dots, a_{j,n}]$ be the  global section type of $\FF_{H_j}$ (as a sheaf on $H_j\simeq \PP^j$);  especially $\st(\FF)=\mathbf{a}_{1}=\mathbf{b}$.
Then
\begin{multline}\label{grossa}
h^0 \OPN(\mathbf{b}) -h^0\FF\geq  \\ 
 \geq \sum_{\substack{1\leq i \leq n \\ 2\leq j\leq N \\ a_{j,i}<0}} \hspace{-4pt} h^0\OPN(a_{j-1,i}) 
+ \hspace{-4pt} \sum_{\substack{1\leq i \leq n \\ 2\leq j \leq k\leq N \\ a_{j,i}\geq 0}} \hspace{-8pt}h^0\OO_{\PP^{N-k}}(a_{j,i}) h^0\OO_{\PP^k}(a_{j-1,i}-a_{j,i}-1)
\end{multline}
\end{theorem}

\begin{proof}
We proceed by induction on $N$ and $b_1$.
If $N=1$, the statement is trivial.
If $b_1<0$ then both sides of the inequality are 0.

\medskip

Now we assume the thesis for torsion free sheaves on  $\PP^{N-1}$ and torsion free sheaves whose first integer of the splitting type is $<b_1$.

We immediately have, from Proposition \ref{utile}, the following inequality
\begin{multline}
h^0 \OPN(\mathbf{b}) -h^0\FF
 \geq(h^0 \OPN(\mathbf{b}-1) -h^0\FF(-1))+\\ 
 +(h^0 \OO_H(\mathbf{b}) -h^0\FF_H)+\sum_{\substack{1\leq i\leq n \\ a_{N,i}<0}}h^0\OO_{\PP^{N-1}}(a_{N-1,i})
\end{multline}

We now conclude using the inductive hypothesis on both the dimension of the projective space (for what concernes  $h^0 \OO_H(\mathbf{b}) -h^0\FF_H$) and  on the first integer $b_1$ of the splitting type  (for what concernes $h^0 \OPN(\mathbf{b}-1) -h^0\FF(-1)$) and reckoning that $h^0\OO_{\PP^{r}}(a-1)+h^0\OO_{\PP^{r-1}}(a)=h^0\OO_{\PP^r}(a)$.

\end{proof}

\begin{corollary}\label{risH0}
Let $\FF$ be a rank $n$ torsion free sheaf on $\PP^N$ generated by global sections, or more generally  such that  $\gsr(\FF)=\rk(\FF)$. If $\mathbf{b}=[b_1, \dots, b_n]$ is its splitting type and    $\mathbf{a}=[a_1, \dots, a_n]$ is the global section type of a general plane section, then: 
\begin{equation}\label{menogrande2} h^0\FF\leq h^0\OO_{\PP^N}(\mathbf{b})- \sum_{i=1}^n h^0\OO_{\PP^{N-2}}(a_i)\  h^0\OO_{\PP^2}(b_i-a_i-1)\end{equation}
 as well
\begin{equation}\label{menogrande} h^0\FF\leq h^0\OO_{\PP^N}(\mathbf{b})-  h^0\OO_{\PP^2}(\mathbf{b}-\mathbf{a}-1)\end{equation}
\end{corollary}

\begin{proof} First of all observe that $\gsr(\FF)=\rk(\FF)$ implies $a_i\geq 0$ for every $i$. Thus the inequality (\ref{menogrande2}) is a straightforward consequence of (\ref{grossa}), where we omit in the right side every summand except those with $j=k=2$;  (\ref{menogrande}) can be obtained by the previous one  because in our hypothesis  $h^0\OO_{\PP^{N-2}}(a_i)\geq 1$. 
\end{proof}

For every torsion free sheaf $\FF$,  all the sheaves $\FF(t)$  fulfill  the hypothesis of the previous result for sufficiently hight values of $t$. If we consider the  inequalities (\ref{menogrande2}) for such sheaves, we  obtain an upper bound for $h^0\FF(t)$ set up   by  $h^0\OPN(\mathbf{b}+t)$, which is a polynomial in $t$ of degree $N$,  minus a \lq \lq smaller\rq\rq \ piece given by a polynomial in $t$  of degree $N-2$. We will see in \S 4 that this is a quite reasonable bound,  since the difference $h^0\OPN(\mathbf{b}+t)- h^0\FF(t)$ is in fact a polynomial of $t$-degree $N-2$.

\begin{corollary}\label{primospezz}
For every  torsion free sheaf  $\FF$ on $\PP^N$ with splitting type $\mathbf{b}=[b_1, \dots, b_n]$ and  global section type $\mathbf{a}=[a_1, \dots, a_n]$ the following are equivalent:
\begin{enumerate}
	\item $\FF\simeq \OPN(\mathbf{b})$
	\item $\mathbf{a}= \mathbf{b}$
	\item $h^0\FF(t)=h^0\OPN(\mathbf{b}+t)$ for every $t\in \ZZ$
	\item $h^0\FF(t)=h^0\OPN(\mathbf{b}+t)$ for some $t\geq -a_n$
\end{enumerate}
\end{corollary}
\begin{proof} \textit{1.} $\Rightarrow$ \textit{2.} and \textit{3.} $\Rightarrow$ \textit{4.} are    obvious.

\medskip

  \textit{2.} $\Rightarrow$ \textit{3.}\  By the definition of  global section type we get $h^0\FF(t)\geq h^0\OPN(\mathbf{a}+t)$ and then, in our hypothesis $h^0\FF(t)\geq h^0\OPN(\mathbf{b}+t)$. On the other hand Theorem \ref{H0} gives the opposite inequality.
  
  \medskip
  
 \textit{3.} $\Rightarrow$ \textit{4.}\  As  $t\geq -a_n$, every integer in the global section type of $\FF(t)$ is  positive; moreover the same holds for   the restriction of $\FF(t)$ to a general linear space of every dimension, thanks to Lemma \ref{lemmafacile} \textit{3.}. Using the same notation as in Theorem \ref{H0},  (\ref{grossa}) can be written in the following shorter form:
 \begin{equation}\label{rigrossa}
h^0 \OPN(\mathbf{b}+t) -h^0\FF(t)\geq  
  \sum_{\substack{1\leq i \leq n \\ 2\leq j \leq k\leq N }} h^0\OO_{\PP^{N-k}}(a_{j,i}+t) h^0\OO_{\PP^k}(a_{j-1,i}-a_{j,i}-1)
\end{equation}
Thus condition \textit{4.} can hold   only if $a_{j-1,i}=a_{j,i} $ for every $i,j$, that is if $\mathbf{a}=\gst(\FF)=\gst(\FF_H)=\dots = \st(\FF)=\mathbf{b}$.
\medskip

\textit{2.+ 3.} $\Rightarrow$ \textit{1.}\  By definition of global section type and the hypothesis $\mathbf{a}=\mathbf{b}$ there is  an injective map $\phi \colon \OPN (\mathbf{b})\hookrightarrow \FF$;  by the assumption \textit{3.} $\coker(\phi) (t)$  has no global sections for every $t\in \ZZ$. Then $\phi$  is  an isomorphism.
\end{proof}

The following result generalizes Lemma 1.3.3 of \cite{okverde}, where only vector bundles generated by global sections are considered: note that the condition \lq\lq generated by global sections \rq\rq\ implies, in our terminology,  $\mathbf{a}\geq 0$. For another proof see also \cite{prec}.

\begin{corollary}\label{louso} Let $\FF$ be a torsion free  sheaf  on $\PP^N$ with $\gst(\FF)=\mathbf{a}=[a_1, \dots,a_n]$. 

Then   $c_1\geq \sum a_i$ and equality holds if and only if $\FF\simeq \OPN(\mathbf{a})$.

Especially, if $\gsr(\FF)=\rk(\FF)$, then $c_1\geq 0$ and equality holds if and only if $\FF\simeq \OPN^n$
\end{corollary}

\begin{proof}Let $\st(\FF)=[b_1, \dots, b_n]$ and $ \gsr(\FF)=[a_1,\dots,a_n]$. Lemma \ref{lemmafacile} \textit{3. } says that $b_i\geq a_i$ for every $i$; then   $c_1=\sum b_i\geq \sum a_i $. Moreover equality can hold if and only if   $\mathbf{a}=\mathbf{b}$ and we conclude by Corollary \ref{primospezz}.\\
With the further hypothesis $\gsr(\FF)=\rk(\FF)$, we immediately have $c_1\geq 0$ and again for equality we conclude by Corollary \ref{primospezz}.
\end{proof}

\section{Horrocks' splitting criterion}

In this section we will apply the same approach developed in \S 2 in order to generalize the Horrocks' splitting criterion to torsion free sheaves. The results of the two sections are closely related  so that a reference to the previous one (mainly to Corollary \ref{primospezz}) could  shorten the proofs in the present; however we prefer to avoid any reference and present here independent and   self-included proofs in order to underline how simple, direct and effective the employed  methods are.

For other generalizations of the Horrocks criterion,  see also \cite{ABE}.

\begin{theorem}\label{previous}   Let  $\FF$ be a torsion free sheaf  on $\PP^N$, $N\geq 2$, $H$ be a hyperplane in $\PP^N$ and $\mathbf{b}=[b_1, \dots,b_n]$ be an ordered sequence of integers.

Then  $\FF\simeq \OPN (\mathbf{b})$ if and only if $\FF_H\simeq \OO_H(\mathbf{b})$ and either one of the following conditions holds:

	1.   $N=2$   and  $H^1 (\FF(t))=0$ for every $t\leq -b_n$
	
	2. $N\geq 3$    and  $H^1 (\FF(t))=0$ for every $t\ll 0$. 
\end{theorem}

\begin{proof} We only prove the non trivial part, that is we assume that either \textit{1.} or \textit{2.} holds and prove that $\FF$ splits. Up a twist, we may assume $b_n=0$, so that $\OO_H(\mathbf{b})$ is generated by global sections.  For every $t\in \ZZ$, the exact sequences:
\begin{equation}\label{restrizione}0 \rightarrow \FF(t-1)\rightarrow \FF(t) \rightarrow \FF_H (t) \simeq \OO_H(\mathbf{b}+t)\rightarrow 0.\end{equation}
Now observe that the hypothesis imply that $H^0\FF(t)\rightarrow H^0\OO_H(\mathbf{b}+t)$ is surjective. In fact, if $N\geq 3$ the finiteness of $H^1_*(\FF)$ and $H^1_*(\FF_H)=0$ imply  $H^1_*(\FF)=0$; if $N=2$,  the surjectivity of the map is a direct consequence of    \textit{1.} if  $t\leq 0$ and so  also if $t>0$ because $\OO_H(\mathbf{b})$ is generated by global sections.  Thus  (\ref{restrizione}) induces in cohomology the exact sequence:
\begin{equation}\label{alposto}0 \rightarrow H^0\FF(t-1)\rightarrow H^0 \FF(t) \rightarrow H^0 \OO_H(\mathbf{b}+t)\rightarrow 0.\end{equation}

 Now we can fix for every $j$ a global section in $H^0\FF(-b_j)$ whose image is the canonical generator of the $j-$th summand of  $\OO_H(\mathbf{b}-b_j)$ and define a map
$
 \psi \colon \OPN(\mathbf{b})\rightarrow \FF
$: we conclude by showing that this map is in fact an isomorphism.

First of all, $\psi$ is injective because its kernel is a torsion free sheaf which vanishes if restricted to $H$. 
Moreover, summing up on $t\leq t_0$ the equalities $h^0\FF(t)-h^0\FF(t-1)=h^0\OO_H(\mathbf{b}+t)$, we get $h^0\FF(t_0)=h^0\OPN(\mathbf{b}+t_0)$ for every $t_0\in \ZZ$ (note that $H^0\FF(t)=0$ if $t\ll 0$, because $\FF$ is torsion free). Then there is an exact sequence:
\[0 \rightarrow \OPN(\mathbf{b})\rightarrow \FF \rightarrow \TT  \rightarrow 0\]
where $\TT$ is the cokernel of $\psi$. By the exact sequence of cohomology, we can see that, for every $t\in \ZZ$, $H^0\TT(t)=0$,  because $h^0\FF(t)=h^0\OPN(\mathbf{b}+t)$ and $h^1\OPN(\mathbf{b}+t)=0$; then   $\TT=0$ that is  $\psi$ is also surjective and we conclude.
 \end{proof}

\begin{corollary}\label{horrockstf}
Let $\FF$ be a torsion free sheaf on $\PP^N$, $N\geq 2$.  Then:
\begin{center}
 $\FF$ is split $\Longleftrightarrow$
 $H^i_*\FF=0$ for $i=1,\dots,N-1$.
\end{center}
\end{corollary}

\begin{proof}   If  $N=2$,  $\FF_H$ is a torsion free sheaf on $H\simeq \PP^1$ and then it is split; so we conclude thanks to Theorem \ref{previous}. If $N>2$, the restriction of $\FF$ to a general hyperplane  is a torsion free sheaf on $H\simeq \PP^{N-1}$ and has trivial intermediate cohomology too. So,  by induction, $\FF_H$ is split and again we conclude thanks to Theorem \ref{previous}. \end{proof}

\begin{corollary} Let $\FF$ be a reflexive sheaf on $\PP^N$, $N\geq 3$. Then:
\begin{center}
 $\FF$ is split $\Longleftrightarrow$ $\exists$ a plane $H$ such that $\FF_H$ is split.
\end{center}
 \end{corollary}
\begin{proof} By induction we can assume that  $\FF_{H}\simeq \OO_H(\mathbf{b})$, where $H$ is a hyperplane; then  $H^i_*\FF_H=H^i_*\OO_H(\mathbf{b})=0$ for $i=1,\dots ,N-2$ and moreover there is  the exact sequence (\ref{restrizione}). Passing in cohomology and bearing in mind that $H^1\FF(t)=0$ for $t\ll 0$ (see Notation \ref{reflexive}))
we can easily deduce that  $H^i_*\FF=0$ for $i=1,\dots ,N-1$.  So  Corollary \ref{horrockstf} leads  to conclude that $\FF\simeq \OPN^n$, against the assumption.
\end{proof}


\section{$c_2(\FF)\geq c_2(\OPN(\mathbf{b}))$}

In the present section we will show that the split bundle $\OPN (\mathbf{b})$ has the smallest   second Chern class among  torsion free sheaves having the same  splitting type  
$\mathbf{b}=[b_1,\dots,b_n]$. More precisely, as $c_2(\OPN(\mathbf{b}))= \sum_{i<j} b_ib_j$, we will show that the second Chern class of   a rank $n$ torsion free sheaf $\FF$ on $\PP^N$ is strictly bigger than $ \sum_{i<j}b_ib_j$ unless $\FF$ is split. The proofs directly lean on the main results of the previous sections, mainly of \S 2: the connection is given by the  Euler characteristic   $\chi (\FF)=\sum_{i=1}^N (-1)^i
h^i(\FF)$. The equality  $\chi (\FF(t))=h^0\FF(t)$, that holds for every $t\gg 0$, allows us to translate results about  $h^0\FF$ in terms of   Chern classes of $\FF$. In fact by Grothendieck-Hirzebruch-Riemann-Roch (see \cite{HAG} Appendix A, 4.1) $\chi(\FF)$ can be expressed as a polynomial in the Chern classes; for instance if $N=2$, the formula is:
\begin{equation}\label{polinomio}
\chi(\FF)=\frac12 (c_1^2+3c_1)-c_2+n
\end{equation}
We underline that though this argument seems to be bound to  \lq\lq hight twists\lq\lq\ of $\FF$ in order to use the above quoted equality $\chi (\FF(t))=h^0\FF(t)$, however the properties of $c_2(\FF)$ we will obtain also hold for every torsion free sheaf $\FF$, because the difference $c_2(\FF(t))-c_2(\OP (\mathbf{b}+t))$ does not depend on $t$.

\begin{theorem}\label{secondobound}
Let $\FF$ be a rank $n$ torsion free sheaf on $\PP^N$, with $\st(\FF)=[b_1,\dots,b_n]$ and let $[a_1,\dots,a_n]$ the  global section type
 of its restriction to a general plane. Then:
\[
c_2(\FF)\geq \sum_{1\leq i<j\leq n}b_ib_j+\sum_{i=1}^n\frac{(b_i-a_i)(b_i-a_i+1)}{2}.
\]
and  \[ c_2(\FF)=\sum b_ib_j \Longleftrightarrow \FF\simeq \OPN(\mathbf{b}).\]
\end{theorem}
\begin{proof}
As $c_2(\FF)=c_2(\FF_{\PP^2})$, we may assume that $\FF$ is defined on $\PP^2$. Using  (\ref{polinomio}) for both sheaves $\FF$ and   $\OP(\mathbf{b})$ that  have the same rank $n$ and first Chern class $c_1=\sum b_i$, we find: \[\chi(\OP(\mathbf{b}))-\chi(\FF)=c_2(\FF)-c_2(\OP(\mathbf{b}))=c_2(\FF)-\sum_{ i<j}b_ib_j\]
Moreover, as both integers
$c_2(\FF)-\sum b_ib_j$ and $\sum \frac{(b_i-a_i)(b_i-a_i+1)}{2}$ are invariant by twist,  we can also assume  that $\chi(\OP(\mathbf{b}))-\chi(\FF)=h^0\OP(\mathbf{b})-h^0\FF$.  So  Corollary \ref{risH0} allows to conclude. 
\end{proof}

\begin{corollary}\label{c2positiva}
Let $\FF$ be a  non-split torsion free sheaf on $\PP^N$, with $\st(\FF)=[b_1,\dots,b_n]$. 
 If either $b_1\leq 0$ or $b_n\geq 0$, then $c_2(\FF)> 0$. 
 So, there are at most $b_1-b_n-1$  integers $t$ such that $c_2(\FF(t))\leq 0$.
\end{corollary}

Now we consider  semistable sheaves. Recall that a torsion free sheaf on $\PP^N$ is semistable  if for every non-zero subsheaf $\GG\subset \FF$ we have $\displaystyle{ \frac{c_1(\GG)}{\rk (\GG)}}\leq \frac{c_1(\FF)}{\rk (\FF)}$; if moreover the inequality is strict for every $\GG$ such that  $\rk (\GG)<\rk (\FF)$ then $\FF$ is stable.

 Now we will prove that  every non-split semistable torsion free sheaf of any rank on $\PP^2$ has a positive second Chern class.

\begin{lemma}\label{poche} Let  $\FF$ be a non-split rank $n$ torsion free sheaf on  $\PP^N$, normalized so that  $1-n\leq c_1(\FF) \leq 0$. If $\FF$ is semistable, then:
\begin{itemize}
\item[1.]  $h^0(\FF(-1))=0$   
\item[2.]  $h^0(\FF)=0$ if $c_1<0$ and  $h^0(\FF) \leq n-1$  if $c_1=0$;
\end{itemize}
\end{lemma}

\begin{proof}
 {\textit{1.}} and the first part of \textit{2.} directly follow from semistability (see also  \cite{okverde} Lemma 1.2.5 and Remark 1.2.6). Moreover the statement obviously holds if $n=1$. 
 Then assume $n\geq 2$,  $c_1=0$ and $h^0\FF \geq n$; we may also suppose that $n$ is the minimal one for which such a sheaf exists. 
 
Take $n$ linearly independent global sections of $\FF$ and define the sheaf  $\HH $  generated by them, so that:
$0 \rightarrow \HH \rightarrow \FF
 $  and $\OPN^n \rightarrow \HH \rightarrow 0 $.
Since   $\FF$ and $\OPN^n$ are  semistable,  we must have both  $c_1(\HH)\leq 0$ and  $c_1(\HH)\geq 0$, so that $c_1 (\HH )=0$; as a consequence also $\HH$ have to be semistable, because $\FF$ is.  Moreover, by construction, $h^0\HH\geq n$. Thus, the minimality of $n$ forces  $\rk(\HH)$ to be $n$, so that the above surjective map must be an isomorphism:   $\HH\simeq \OPN ^n$.
So there is an injective map $\OPN^n \hookrightarrow \FF$ and Corollary \ref{louso} allows to conclude.
\end{proof}

\begin{lemma}\label{duality} Let  $\FF$ be a  torsion free sheaf on  $\PP^2$. Then $h^2\FF-h^1\FF \leq h^0\FF^{\vee}(-3)-h^1\FF^{\vee}(-3)$.
\end{lemma}
\begin{proof} 
Since $\FF$ is torsion free, the canonical morphism $\FF\rightarrow \FF^{\vee\vee}$ is injective and its cokernel is supported by a finite subset of $\PP^2$. Moreover, $\FF^\vee$ is locally free.

Applying Serre duality to $\FF^\vee$ one gets that $h^2\FF=h^2\FF^{\vee \vee}=h^0(\FF^\vee(-3))$ and $h^1\FF\geq h^1\FF^{\vee\vee}=h^1(\FF^\vee(-3))$.
\end{proof}

\section{Applications to the discriminant}

For every rank $n$ torsion free sheaf $\FF$ on $\PP^N$  the discriminant is the number $$ \Delta(\FF)=2nc_2-(n-1)c_1^2 $$ which is  invariant both up twists and dual, that is $\Delta(\FF)=\Delta(\FF(t))=\Delta(\FF^{\vee})$. 

\begin{theorem}\label{c2stabile} Let $\FF$ be a non split semistable torsion free sheaf on $\PP^2$. Then:
\[\Delta (\FF) \geq 2n.\]
Furthermore, if $\FF$ is stable then:
\[\Delta (\FF) \geq  \frac34 n^2.\]
\end{theorem}

\begin{proof}  We may assume  that  $-\frac{n}{2} \leq c_1(\FF) \leq 0$,  substituting to $\FF$, if necessary, a suitable  twist or   dual of a twist. By Lemma \ref{duality}  and Lemma \ref{poche} we have $ h^2\FF-h^1\FF \leq h^0\FF^{\vee}(-3)= 0$; moreover $h^0(\FF) =0$ if $\FF$ is stable and $h^0(\FF) \leq n-1 $ if $c_1=0$ and $\FF$ is properly semistable. Then the Euler formula (\ref{polinomio}) gives:
\begin{description}
\item[$1.$] $c_2\geq  \frac{(c_1^2+3c_1)}{2}+n $ if $\FF$ is stable 
\item[$2.$] $c_2 \geq 1 $ if $\FF$ is properly semistable.
\end{description}

If $\FF$ is stable $1.$  gives $\Delta (\FF) \geq  (c_1+n)(c_1+2n)$ and the right hand assumes its minimum when $c_1$ assumes its minimum $-\frac{n}{2}$.
 
 If $\FF$ is properly semistable $2.$   gives directly the wanted relation $\Delta (\FF) \geq  2n$ 
\end{proof}

 The positivity of the discriminant   for semistable torsion free sheaves on a smooth surface was proved by Bogomolov; for a proof see \cite{HL}, Theorem 3.4.1.
 
\begin{corollary} If a plane section $\FF_H$ of a non totally split torsion free sheaf $\FF$ is semistable, then $c_2(\FF) > 0$.
\end{corollary}
 For rank $2$ vector bundles on $\PP^2$, the  above result is a consequence of Schwarzenberger inequality $c_1^2-4c_2 \leq 0$; it can also be extended to rank $2$ reflexive sheaves on $\PP^N$ as their general plane section is a semistable vector bundle (see \cite{okverde} or \cite{stable}).  If $\rk (\FF) >2$ it is a difficult open question to state if the general plane section of a semistable vector bundle is semistable too.  What it is known is that the splitting type of a semistable vector bundle on $\PP^N$ has no gap (\cite{okverde} Theorem 2.1.4 and Corollary 1).
 Recall that an ordered sequence of integers $\mathbf{b}=[b_1, \dots, b_n]$ has no gap if $b_i-b_{i+1}\leq 1$ for every $i=1, \dots, n-1$.  Using this property of the splitting type and Theorem \ref{secondobound} we can obtain a generalization of  Schwarzenberger inequality to semistable bundles of   rank $n$ on $\PP^N$ (see \cite{okverde} Problem 1.4.1). The proof is just a  numerical computation; the following example shows as the general proof works.

\begin{example} Let $\FF$ be a rank $n \leq 4$ semistable vector bundle  (or, more generally,  torsion free sheaf   whose splitting type has no gap) on $\PP^N$.
If $\FF$ is not totally split, then $c_2\geq 0$.

 Let in fact  $\mathbf{b}=[b_1, \dots, b_n]$ be the splitting type of $\FF$. Theorem \ref{secondobound} says that $c_2(\FF)$ can be negative only if $b_1 > 0 >b_n$. So it will be sufficient to check every possible splitting type with a positive $b_1$ and a negative $b_n$ and use Theorem \ref{secondobound}.  
\begin{itemize}
	\item If $n=2$, no such case can realizes.
	\item if $n=3$ the only possible splitting type is $[1,0,-1]$. Then $c_2 > \sum b_ib_j=-1$.
	\item If $n=4$ there are 5 possible cases: three of the integers $b_i$ are  $1,\ 0,\ -1$ while the forth can be any element $b\in \{ -2,\ -1, \ 0,\ 1,\ 2\}$. Then $c_2>\sum b_ib_j=-1+b(-1+0+1)= -1$.\end{itemize}
\end{example}

\begin{theorem}\label{compleanno} Let $\FF$ be a rank $n$ semistable vector bundle on $\PP^N$ or, more generally, a rank $n$ torsion free sheaf on $\PP^N$  whose splitting type has no gap. If $c_i=c_i(\FF)$, then:
\begin{equation}\label{mattina} \Delta (\FF) \geq -\frac{1}{12}n^2(n^2-1).\end{equation}
More precisely, if $\overline{c}$ is the only integer such that $0\leq \overline{c}\leq \frac{n}{2}$ and $c_1\equiv \pm \overline{c} \bmod n$, then:
\begin{equation}\label{domenica} \Delta (\FF)  \geq \left\{  \begin{array}{ll} \displaystyle{-\frac{2n}{4} {n\choose 3}}- (n-1)\overline{c}^2  & \hbox{ if $n$ is even} \\ \\
\displaystyle{-\frac{2n}{4}{n+1\choose 3}} +(n-1)\overline{c}(n-\overline{c}) & \hbox{ if $n$ is odd}\end{array}\right. \end{equation}
Moreover equality can hold in either inequality only if $\FF$ is totally split;  especially in the  semistable case,   only if $\FF\simeq \OPN$.
\end{theorem}

\begin{proof} The left hand of the three inequalities is invariant by twist and dual; so it is sufficient to prove the statement assuming that $c_1=\overline{c}$. Thus (\ref{domenica}) becomes:
\begin{equation}\label{domenica1} c_2  \geq \left\{  \begin{array}{ll} \displaystyle{-\frac14 {n\choose 3}} & \hbox{ if $n$ is even} \\ \\
\displaystyle{-\frac14 {n+1\choose 3}} +\frac12 (n-1)\overline{c} & \hbox{ if $n$ is odd}\end{array}\right. \end{equation}

By hypothesis, the splitting type $\mathbf{b}=[b_1, \dots, b_n]$ of $\FF$ has no gap (for vector bundles see see \cite{okverde} Theorem 2.1.4 and Corollary 1). Moreover we know that $c_2\geq \sum b_ib_j=\frac12 (\sum b_i)^2 -\frac12 \sum b_i^2=\frac12 \overline{c}^2 -\frac12 \sum b_i^2$ (Theorem \ref{secondobound}). 

Now let us  compute the maximal value of $\sum b_i^2$ for any sequence $[b_1, \dots, b_n]$ without gaps and such that $\sum b_i=\overline{c}$. 

Clearly the maximal value is reached when
$\mathbf{b}$ is  either the sequence $[-\frac{n-2}{2}, \dots, \frac{n-2}{2}]$ plus one more item given by $\overline{c}$ if $n$ is even, or  the sequence $[-\frac{n-3}{2}, \dots, \frac{n-1}{2}]$ plus one more item given by $\overline{c}-\frac{n-1}{2}$ if $n$ is odd.

The formula for the sum of consecutive squares:
\[\sum_{i=1}^k i^2=\frac{2k^3+3k^2+k}{6}\]
and a straightforward computation give (\ref{domenica1}). 

In order we deduce (\ref{mattina}) in the odd case, we only forget the last summand, while in the even one we observe that the maximal value of $-\frac{(n-1)}{2n}\overline{c}^2 $ is obtained when $\overline{c}=\frac{n}{2}$. 

Finally, equality can hold only if $\FF\simeq \OPN(\mathbf{b})$ and  $\mathbf{b}$ is either one of the two special splitting type above considered; the only semistable  split bundle  of that  type is $\OPN$. \end{proof}

For a similar result for semistable torsion free sheaves see \cite{HL} Theorem 7.3.1.

\section{$c_s(\FF(t))$ and $c_s(\OPN(\mathbf{b}+t))$}

Let   $\FF$  be a rank $n$ torsion free sheaf on  $\PP^N$ with Chern classes $c_i$. 

An argument similar to the one of $c_2$  can be applied to each Chern class $c_s$, $s\leq n$. Using (\ref{cherncontwist}) we can see that   $c_s(\FF(t))$ is given by  a degree   $s$ polynomial whose leading coefficient only depends on  $n$ and $s$ and whose next coefficient only depends on $n$, $s$ and $c_1$. So, if  $\mathbf{b}=\st(\FF)$, the difference  $c_s(\FF(t))-c_s(\OPN(\mathbf{b}+t))$ is a polynomial  of degree lower or equal to $s-2$, because the two sheaves have both the  same rank and and the same first Chern class. Moreover the coefficient of $t^{s-2}$ is $\left({n-2}\atop{s-2}\right)(c_2(\FF)-\sum b_ib_j)$, that cannot vanish, in fact it is strictly positive,  unless $\FF$ is split (see Theorem \ref{secondobound}). We can summarize the consequences in the following

\begin{corollary}\label{alte} Let $\FF$ be a rank $n$   torsion free sheaf on $\PP^N$ with splitting type  $\mathbf{b}$. 

 If    $3\leq s \leq min\{ n,\ N \}$, then $c_s(\FF(t)) > c_s(\OPN(\mathbf{b}+t))$ for $t\gg 0$. More precisely:  \[\Lambda^{(s)}_{\FF}(t):= c_s(\FF(t))-  c_s(\OPN(\mathbf{b}+t)) \] is  a degree $s-2$ polynomial with a positive leading coefficient:
\begin{equation}\label{ultimau} \Lambda^{(s)}_{\FF}(t)=\left({n-2}\atop{s-2}\right)(c_2(\FF)-\sum_{i\neq j} b_ib_j)\ t^{s-2}+ (\hbox{lower terms}).\end{equation}

\end{corollary}

 As far as Chern classes of  codimension $s>n$ are concerned, (\ref{cherncontwist})   shows that $c_s(\FF(t))$ is a polynomial in $t$ of degree   $\leq s-n-1$ and that the coefficient of $t^{s-n-1}$ is $(-1)^{s-n-1} c_{n+1}$. So, if $c_{n+1}=0$, also $c_s=0$ for every $s>n$: it is the case  of vector bundles. On the other hand if $c_{n+1}\neq 0$, then $c_s(\FF(t))$ are definitely positive or definitely negative,  depending on the sign of $c_{n+1}$ and the parity of $s-n-1$.

 \section{Bounds for cohomology and regularity} 
 
  In  the paper \cite{EF} it is proved that  the absolute value of the $i$-th Chern class of every rank $n$ vector bundle on $\PP^N$ is upper-bounded by a polynomial function depending on  $n$,  $N$, the splitting type $\mathbf b$ and the first and second Chern classes $c_1$ and $c_2$.

 Actually, no such a bounding formula exists for the wider class of torsion free sheaves, as we can see in the following example. 
 
\begin{example}  Let $\II_Y$ be the ideal sheaf of any set of $M$ different points in $\PP^N$, which is of course a torsion free sheaf of rank 1 and generic splitting type $[b_1=0]$. An easy computation shows that its Chern classes vanish, except  the last one $c_N(\II_Y)=\pm M(N-1)!$. 

For every integer $n$  consider the rank $n$ torsion free sheaf $\FF=\II_Y \oplus \OPN^{n-1}$. Then $\st(\FF)=\mathbf{b}=[0, \dots,0]$ and $c_i(\FF)=0$ for every $i=1, \dots, N-1$, while $ c_N(\FF) =\pm M(N-1)!$.  It is evident that $\vert c_N(\FF)\vert$ cannot be bounded by a function (either polynomial or not) $f(n,\ N,\  c_1,\ c_2,\ \mathbf{b})$ because,   for fixed $n$ and $N$, $f(n,\ N,\  0,\ 0,\ \mathbf{0})$ is a constant, while $\vert c_N(\FF)\vert $ can be as big as we want.
\end{example}

However we can generalize   Theorem 3.3 and Theorem 4.2 of \cite{EF} to reflexive sheaves; our proof follows in its essential points the proof for  vector bundles.

In the following $\FF$ will be a rank $n$ reflexive sheaf on $\PP^N$, $\delta$   the number $\delta_2=c_2-\sum b_ib_j$  and $d$ the {\it diameter}  $d=b_1-b_n$  of the splitting type $\mathbf b$; both $\delta$ and $d$ are invariant up any twist of $\FF$ .

\begin{lemma}\label{sottoinsieme} 
Fixing $n, N$ non-negative integers, every polynomial function depending on $\{ c_1,c_2, \mathbf{b}, d, \delta_2\}$ can be  bounded by polynomials functions depending on either one  of the following sets:  $\{ \mathbf{b},  \delta_2\}$, $\{ \mathbf{b},  c_2\}$,
 $\{ c_1,c_2,  d\}$,
$\{ c_1, d, \delta_2\}$. 
\end{lemma}

\begin{proof}
It is sufficient to recall that $c_1=\sum b_i$, $\delta_2= c_2-\frac12\left(  \sum b_i\right)^2 +\frac12 \sum b_i^2$ and observe that the following inequalities hold:
$$\frac{c_1}{n}-d\leq b_n \leq \dots \leq b_1 \qquad \text{and} \qquad \frac{c_1}{n}+d\geq b_1 \geq \dots \geq b_n.$$
\end{proof}

 \begin{theorem} \label{polinomi}
 For any choice of non-negative integers $n$, $N$ and $s$, with $N\geq 2$ and $3\leq s \leq N$, there are polynomial functions $P_{n,N}$, $Q_{n,N}$, $C_{n,N,s}$  depending on $\{ c_1,c_2, \mathbf{b}, d, \delta_2\}$ (or on one of the sets of Lemma \ref{sottoinsieme}), such that for any reflexive sheaf $\FF$ of rank $n$ on $\PP^N$ we have:
 \begin{enumerate}
	\item[1.] $h^i\FF \leq P_{n,N}$ for all $0\leq i\leq N$;
	\item[2.]  $h^i\FF(k)=h^1\FF(-k)=h^0\FF(-k)=0$ for all   $k \geq Q_{n,N}$ and $1\leq i\leq N$;
	\item[3.] the Castelnuovo-Mumford regularity of $\FF$ is lower than  $ Q_{n,N}$ and especially $\FF(k)$ is generated by global sections for $k \geq Q_{n,N}+N$;
	\item[4.] $\vert c_s(\FF)\vert  \leq C_{n,N,s}$.
\end{enumerate}
 \end{theorem}

\begin{proof} 

First of all, we observe that \textit{3.} is a consequence of \textit{2.} for Castelnuovo-Mumford (\cite{M}, Proposition on page 99). \\
Furthermore we can obtain \textit{4.} from \textit{1.}proceeding by induction on $s$ and $N$; in fact, restricting the sheaf to a generic linear space of dimension $s$, we can apply Riemann-Roch for sheaves on $\PP^s$ and write $c_s$ as a polynomial in the Chern classes  $c_1, \dots, c_{s-1}$ and the Euler characteristic $\chi (\FF)$.

Point \textit{1.} for $i=0$ and $i=N$ and point \textit{ 2.} for $i=N$ follow immediately from (\ref{menogrande}). 

We have to prove \textit{1.} and \textit{2.} for $1\leq i\leq N-1$.  We  start with $N=2$.

In this case, $h^1\FF=h^0\FF+h^2\FF-\chi (\FF)$  and so $h^1\FF$ is bounded by a polynomial depending on one of the sets of invariants considered above.

In order to prove \textit{2.} for $i=1$ we consider a general line $L$ and we observe that for all $k\geq  -b_n$,     $\FF_L(k)$ is generated by global sections and $h^1\FF_L(k)=0$. If there is $k_1\geq -b_n$  such that the map  $H^1\FF(k_1-1) \to H^1\FF(k_1)$ is not only surjective but also injective, then it is an isomorphism for all $k\geq k_1$. Then, for $k\geq -b_n$, the sequence $h^1\FF(k)$ is strictly decreasing or definitely constant. Because for Serre $h^1\FF(k)=0$ for $k\gg 0$, then $h^1\FF(k)=0$ for all $k \geq -b_n+h^1\FF $.\\
Finally, since $\FF$ on $\PP^2$ is a vector bundle, using duality we have  $h^1\FF(-k)=h^1\FF^{\vee}(-k-3)$ and so the last condition holds applying what already shown to the sheaf $\FF^\vee$, whose invariants agree with those of $\FF$ up to the sign.

\medskip

Consider now $N\geq 3$. We prove \textit{1.} and \textit{2.} for  $1\leq	 i\leq N-1$ proceeding by induction on $N$ and assuming the thesis for a generic hyperplane $H\cong\PP^{N-1}$.\\

For $i>1$, we have $h^i\FF(k) = h^i\FF(k+1)$ for all $k\geq k_0=Q_{n,N-1}$ and $ h^i\FF(k)=0$ for $k\gg 0$; then $h^i\FF(k)=0$ for all $k\geq k_0$. Furthermore $h^i\FF(k)\leq h^i\FF(k+1)+h^{i-1}\FF_H(k+1)$; we then obtain a bound for $h^i\FF$ summing up the bounds for $h^{i-1}\FF_H(k+1)$, $k$ from 0 to $k_0$.
 
 \medskip
 
 The only case left is $i=1$. We recall that the first cohomology module of a reflexive sheaf is of finite type; indeed, for $k\gg 0$, $n-1$ general sections of $\FF(k)$ degenerate on a codimension 2 subvariety $Y$ given by the exact sequence
$$
0\to\OO_{\PP^N}^{n-1}\to \FF(k)\to \IY(c_1)\to  0 
$$
 (see \cite{okted}, \S 2). $Y$ is generically smooth and without embedded or isolated components of codimension $\geq 3$. 
 
 Using the same  arguments as before, we prove that $h^1\FF(-k)=0$ for all $k \geq k'$ such that $h^0\FF_H(-k)=h^1\FF_H(-k)=0$. Using $h^1\FF(-k)\leq h^1\FF(-k-1)+h^1\FF(-k)$, we obtain a bound for $h^1\FF$ summing up the bounds for $h^1\FF(-k)$, $k$ from 0 to $k'$.
\\
Finally, the sequence of $h^1\FF(k)$'s is strictly decreasing or definitely 0; so $h^1\FF(k)=0$ if  $k\geq k''+ h^1\FF(k'')$.
 \end{proof}

\begin{corollary}

For any choice of non-negative integers $n$, $N$ and $s$, with $N\geq 2$ and $3\leq s \leq N$, there are polynomial functions $P_{n,N}$ and $C_{n,N,s}$  depending on $\{ c_1,c_2, \mathbf{b}, d, \delta_2\}$ (or on one of the sets of Lemma \ref{sottoinsieme}) and $t$, such that for any twist $\FF(t)$ of a reflexive sheaf $\FF$ of rank $n$ on $\PP^N$ we have:
 \begin{enumerate}
	\item[1.] $h^i\FF(t) \leq P_{n,N}$ for all $0\leq i\leq N$;
	\item[2.] $\vert c_s(\FF(t))\vert  \leq C_{n,N,s}$.
\end{enumerate}
\end{corollary}
 
\section{Sharp examples}

In this last section we examine the main results of the paper, in which inequalities are  presented, in order to check how they are sharp.   First of all we take into consideration the bound on the dimension of the  space of global sections given by (\ref{grossa}), which is the starting point for a big amount  of our results. 

\medskip

Let   $\FF$  be a rank $n$ non totally split torsion free sheaf on  $\PP^N$ with Chern classes $c_i$. Recall that the Euler characteristic $\chi (\FF(t))$, as a function of $t$ and $c_i$,  is a polynomial with rational coefficients which depend on $n$ and $N$. By the splitting principle, this polynomial can be computed using free sheaves; in this way we can see that, if $c_i$ is regarded as a variable of weight $i$, $\chi (\FF(t))$ has \lq\lq global degree\rq\rq\  $N$. Especially:
\begin{itemize}
	\item the coefficients of $t^N$ and   $t^{N-1}$ only depend on $n$,  $N$  and $c_1$
	\item the coefficient of $t^{N-2}$ is given by ${-\frac{c_2}{(N-2)!}+} $ \textit{terms  depending on} $ n, N, c_1$.
\end{itemize}

If  $\mathbf{b}$ is the splitting type of $\FF$, we can  measure how much $\FF$ differs from the free sheaf $\OPN(\mathbf{b})$ (in a natural way the closest free sheaf to $\FF$) through the polynomial  $D_{\FF}(t)= \chi(\OPN(\mathbf{b}+t))-\chi(\FF(t))$. As    $\OPN(\mathbf{b})$ has both the same rank and first Chern class as $\FF$,  the degree of $D_{\FF}(t)$   as a polynomial in $t$ can be at most $N-2$; but Theorem \ref{secondobound}   insures that the coefficient of $t^{N-2}$ cannot vanish,   so that the degree of $D_{\FF}(t)$ with respect to $t$ is always $N-2$.
If $t\gg 0$ we have  $\chi(\FF(t))=h^0\FF(t)$ and $\chi(\OPN(\mathbf{b}+t))=h^0 \OPN(\mathbf{b}+t)$, so that also $h^0 \OPN(\mathbf{b}+t)- h^0\FF(t)=D_{\FF}(t)$,   a polynomial of degree $N-2$ with respect to $t$. 

\medskip

On the other hand we can compute the lower bound for $h^0 \OPN(\mathbf{b}+t)- h^0\FF(t)$ given by the inequality  (\ref{grossa}). If $t\gg 0$, the global section type of $\FF(t)$ and those of its restrictions to general linear spaces are all positive, so that we only have to consider the last sum in (\ref{grossa}); moreover thanks to Corollary \ref{primospezz} we know that there is at least one   non-zero contribution for $j=k=2$. Then  the bound given by (\ref{grossa}) is a polynomial of degree $N-2$ with respect to $t$, the same degree as $D_{\FF}(t)$.

Furthermore,  equality can in fact realize in (\ref{grossa}), as shown by the following examples.

\begin{example} Let $\FF$ be a normalized null-correlation bundle on $\PP^3$; it has splitting type  $[0,0]$ and it is linked to the disjoint union of a pair of line $Y$ by the exact sequence:
	\begin{equation}0\rightarrow \OO_{\PP^3}(-1) \rightarrow \FF \rightarrow \II_Y(1) \rightarrow 0
\end{equation} so that $h^0(\OO_{\PP^3}(t)\oplus \OO_{\PP^3}(t))-h^0\FF(t)=t+2$ for every $t\geq 1$. Moreover  $\gst(\FF)=[-1,-1]$ and  $\gst(\FF_H)=[0,-1]$ for a general plane $H$, so that $t+2$ is exactly the bound given by (\ref{grossa}). 
\end{example}

\begin{example}  Let $n$ be any positive integer and let $\mathbf{b'}=[0, \dots, 0]$ be a sequence of length $n-1$. Fix a line $Y$ in $\PP^3$;  a general global section  of $\omega_Y\otimes \OO_{\PP^3}(\mathbf{b'}+4)\simeq \OO_Y(\mathbf{b'}+2)$ defines, as an extension, a rank $n$ reflexive sheaf $\FF$ (see \cite{stable}, Theorem 4.1):
	\begin{equation}\label{ultima2}0\rightarrow \OO_{\PP^3}^n \rightarrow \FF \rightarrow \II_Y\rightarrow 0.
\end{equation}
   It is easy to check that $\gst(\FF)=\gst(\FF_H)=[0, \dots, 0,-1]$ ($H$ a general plane) while $\st(\FF)=\mathbf{b}=[0,\dots,0]$, so that that  (\ref{grossa}) gives the bound $h^0 \OO_{\PP^3}(\mathbf{b}+t)- h^0\FF(t)\leq t+1$ for every $t\gg 0$; using  (\ref{ultima2}) we can see that $h^0 \OO_{\PP^3}(\mathbf{b}+t)- h^0\FF(t)=h^0\OO_{\PP^3}(t)-h^0\II_Y(t)=h^0\OO_{\PP^1}(t)=t+1$  that is  equality holds for every $t\geq 1$. 
   
   Furthermore we can obtain equality for every positive value $k$ of the leading coefficient of $D_{\GG}(t)$,  taking the sheaf $\GG$  direct sum of $k$ copies of the above considered  sheaf $\FF$: in fact $\st(\GG)=[0,\dots, 0]$, $\gst(\GG)=\gst(\GG_H)=[0, \dots, 0,-1,\dots ,-1]$,    $D_{\GG}(t)= kt+k$ for every $t\geq 1$, realizing equality in (\ref{grossa}) for every $t\geq 1$. 
\end{example}

Now we consider the bounds on  $c_2$ and $S_{12}$ given in Theorem \ref{secondobound}  and Theorem \ref{c2stabile}. Observe that the second Chern classes of the sheaves  considered in the above examples  assume the  minimum value allowed for non totally split sheaves, that is $\sum b_i b_j +1$. Now we will see that, on the other hand, $c_2$ can assume every value above that; more precisely for every $n\geq 2$, sequence  $\mathbf{b}=[b_1, \dots, b_{n}]$  and  positive integer $r$,  there is a rank $n$ reflexive sheaf $\FF$ such that $\st(\FF)=\mathbf{b}$ and $c_2(\FF)=\sum b_i b_j +r$.

If we consider a sequence $\mathbf{b}$ corresponding to the minimal value of $S_{12}$ for a given $\sum b_i$, the reflexive sheaves  obtained in Example \ref{estensione} are extremal with respect to the bounds given in Theorem \ref{estensione}.

 As $c_2(\FF)-c_2(\OO_{\PP^N}(\mathbf{b}))$ and $S_{12}(\FF)$ are invariant by twist, we can assume, without lost in generality,   $b_n=0$.

\begin{example}\label{estensione} Let $r$ be any positive integer and $\mathbf{b'}=[b_1, \dots, b_{n-1}]$ be any sequence of non negative integers. Fix a curve $Y$ in $\PP^3$  which is a complete intersection $(r,1)$. A general global section  of $\omega_Y\otimes \OO_{\PP^3}(\mathbf{b'}+4)\simeq \OO_Y(\mathbf{b'}+r+1)$ defines, as an extension, a rank $n$ reflexive sheaf $\FF$ :
	\begin{equation}\label{ultima}0\rightarrow \OO_{\PP^3}(\mathbf{b'}) \rightarrow \FF \rightarrow \II_Y\rightarrow 0.
\end{equation}
Using multiplicativity of Chern polynomials  we can compute Chern classes: $c_2(\FF)=c_2(\OO_{\PP^3}(\mathbf{b'})) +deg(Y)=c_2(\OO_{\PP^3}(\mathbf{b}))+r$, where $\mathbf{b}=[b_1, \dots, b_{n-1},0]$ is the splitting type of $\FF$. 
\end{example}

 The  second item of Corollary \ref{c2positiva}  allows a negative $c_2$ only for $b_1-b_n-1$ twists of a sheaf $\FF$ with splitting type $[b_1, \dots, b_n]$. In the following example are presented sheaves $\FF$ such that $c_2(\FF(t))<0$ if and only if $b_n < t < b_1$ that is for $t$ in the wider interval allowed by  the above quoted result. 

\begin{example} Fix a line $Y$ in $\PP^3$ and  any positive integer $b$. A general global section of $\omega_Y(2+2b)$ defines, as an extension, a rank 2 reflexive sheaf $\FF$:
	\[0\rightarrow \OO_{\PP^3}(b) \rightarrow \FF \rightarrow \II_Y(-b)\rightarrow 0.
\]
The splitting type of $\FF$ is $\mathbf{b}=[b,-b]$  and the first two  Chern classes are $c_1(\FF)=0$ and $c_2(\FF)=-b^2+1$.  Then $c_2(\FF(t))=t^2-b^2+1$  is negative if and only if $-b+1\leq t \leq b-1$.

Note that  for every  $t\geq b$  we have    $h^0\OO_{\PP^3}(\mathbf{b}+t)-h^0\FF(t)=t-b+1$ which is the minimum allowed by Theorem \ref{secondobound};  as a consequence, also  $c_2(\FF)(t)=t^2-b^2+1$ is   the minimum allowed by Corollary \ref{c2positiva}.

In a easy way we can obtain analogous sharp examples of any rank $n$ taking into consideration  for instance the sheaf $\FF\oplus \OO_{\PP^3}^{n-2}$, where $\FF$ is that considered above. 
 \end{example}

\bigskip

\noindent   Cristina Bertone and Margherita Roggero\\
   Dipartimento di Matematica dell'Universit\`{a} \\ 
           Via Carlo Alberto 10 \\ 
           10123 Torino, Italy \\
   {\small cristina.bertone@unito.it} \\ {\small margherita.roggero@unito.it 
}
\end{document}